\documentclass[
	a4paper, %
	8pt, %
]{LTJournalArticle}

\usepackage{etoolbox}
\providetoggle{preprint}
\settoggle{preprint}{true}

\usepackage{authblk}
\newcommand\blfootnote[1]{%
  \begingroup
  \renewcommand\thefootnote{}\footnote{#1}%
  \addtocounter{footnote}{-1}%
  \endgroup
}

\usepackage{siunitx}
\usepackage{amsmath,amsfonts,amssymb}
\usepackage{cleveref}
\usepackage{verbatim}

\usepackage[disable]{todonotes}
\setuptodonotes{inline}

\newcommand{\TheTitle}{
  Optimal dynamic thermal plant control: A study and benchmark
  }

\newcommand{\Secref}[1]{Sec.~\ref{#1}}

\newcommand{\vu}{\ensuremath{\mathbf{u}}}
\newcommand{\vv}{\ensuremath{\mathbf{v}}}
\newcommand{\vy}{\ensuremath{\mathbf{y}}}

\newcommand{\st}{\quad \ensuremath{\text{s.t. \,}}}
\newcommand{\divop}{\operatorname{div}}

\newcommand{\R}{\ensuremath{\mathbb{R}}}
\newcommand{\Uad}{\ensuremath{U_{ad}}}
\newcommand{\Yad}{\ensuremath{Y_{ad}}}

\newcommand{\abs}[1]{\left\lvert#1\right\rvert}

\newcommand{\vcp}{\ensuremath{V_{c/p}}}
\newcommand{\nrn}{\ensuremath{N_V}}

\newcommand{\ecp}{\ensuremath{E_{c/p}}}

\newcommand{\nre}{\ensuremath{N_E}}

\iftoggle{preprint}{
    \title{\TheTitle{}} %
    \author[1]{Thomas Grandits}
    \author[2]{Stefano Coss}
    \author[1]{Gundolf Haase}
    \affil[1]{Department of Mathematics and Scientific Computing, University of Graz, Austria}
    \affil[2]{Arteria Technologies, Austria}
    \date{}
}{}

\begin{document}
\maketitle

\begin{abstract}
    District heating networks play a vital role in thermal energy supply in many countries.
Thus, it comes to no surprise that these has been a central role in improving energy efficiency for private and public energy suppliers alike around the globe.
Many studies have previously investigated the potential of energy saving by low temperature operation of the DHN and the integration of renewable energies.
Many other studies consider this problem in terms of mixed integer linear programming. %
Here, we instead investigate the utilization of well-established continuous optimization methods to improve DHN operation efficiency.
We demonstrate that optimal control is able to model low temperature operation of a DHN for savings of around 8\%, but can even further improve its operation when considering dynamic energy pricing, reducing the cost of operation by roughly 12\%.
We demonstrate the applicability of this method in a realistic, openly available network in Switzerland (OpenDHN), with a total runtime of less than 5 minutes on a standard desktop computer per experiment.
\iftoggle{preprint}{}{
Furthermore, the data and code to reproduce this study will be made publicly available upon acceptance.
}

\end{abstract}

\iftoggle{preprint}{
  \blfootnote{\emph{This work was supported by the Austrian Research Promotion Agency (FFG) over the project \emph{AI4GreenHeatingGrids}  [grant number 923450].}}
}{}
\listoftodos

\section{Introduction}
Many countries worldwide rely on hot water supply through district heating systems for building heating.
While these are currently primarily supplied by fossil fuels, great hopes are put on integrating renewable energy sources into DHNs~\cite{lund_role_2010}.
Modern DHN systems can also be envisioned as low-cost energy storage systems, widely available in many countries~\cite{lund_energy_2016}.
This is especially true for modern 4th generation district heating~\cite{sorknaes_benefits_2020}.
In such systems, one might also be interested in incorporating cheap waste heat sources, or excessively available electrical energy (power-to-heat)~\cite{schweiger_potential_2017}.

To utilize this potential, optimal ways of regulating the heat flow in the network are needed.
In many works, these problems are optimized using mixed integer linear programming optimization methods~\cite{mohring_district_2021,guelpa_optimal_2016,baviere_optimal_2018,krug_nonlinear_2021}.
Here, we want to employ the method of optimal control, which is widely used in engineering disciplines for controlling physical systems governed by ordinary or partial differential equations~\cite{manzoni_optimal_2021}.
At its essence, optimal control deals with finding optimal solutions to a physically constrained problem, which further can be restricted by additional constraints.
It distinguishes between control variables, that influence the physical system, and the state variables, which are physical quantities, influenced by the control variables.
In the context of DHNs, this means we want to optimally control the input temperatures of the heating plants, such that the energy losses in the network are minimized.
Another scenario could be to minimize the energy costs of the network, given a dynamic energy price which is dependent on the time of day.

We will formulate the thermodynamic optimal control problem for a DHN, and how to solve it using modern optimization techniques.
For this, we minimize the thermal ambient losses by imposing the thermodynamic DHN model which will lead to low temperature operation of the DHN.
The problem is further constrained by requiring a minimal heat exchanger input and output temperature at the consumer sites and minimal return temperature at the heating plants with and without dynamic energy pricing.
We will show that this method is able to provide smooth input temperature control while ensuring the constraints imposed on the system.

\section{Methodology}

\begin{table}[htb]
    \centering
    \caption{Description of constants and symbols (first) used throughout this paper along with their description (second) and respective unit and domain (third and fourth).}
    \label{tab:variables}
    \begin{tabular}{cp{3.5cm}cc}
        \hline 
        Sym. & Description & Unit & Domain \\
        \hline \hline
        $c_p$ & Specific heat capacity of water & \unit{\joule \per \kilogram \per \degreeCelsius} & - \\
        $c_A$ & Cross-sectional pipe area & \unit{\meter \squared} & $E$ \\
        $c_d$ & Cross-sectional pipe diameter & \unit{\meter} & $E$\\
        $y$ & Temperature/State & \unit{\degreeCelsius} & $V$ \\
        $y_c$ & Temperature consumed & \unit{\degreeCelsius} & $\ecp{}$ \\
        $y_a$ & Ambient temperature & \unit{\degreeCelsius} & $V$ \\
        $y_p$ & Plant temperature & \unit{\degreeCelsius} & $\vcp{}$ \\
        $t$ & Time & \unit{\second} & - \\
        $\Delta_t$ & Time step & \unit{\second} & - \\
        $\vv$ & Fluid velocity & \unit{\meter \per \second} & $E$ \\
        $\dot{m}$ & Mass flow $\dot{m} = \vv \frac{c_d^2}{4} \pi \rho $ & \unit{\kilogram \per \second} & $E$ \\
        $k$ & Heat transfer coefficient & \unit{\watt \per \meter \per \degreeCelsius} & $E$ \\
        $l$ & Length of pipe segment & \unit{\meter} & $E$ \\
        $\lambda_p$ & Constraint penalty weight & - & - \\
        $\rho$ & Density of water & \unit{\kilogram \per \cubic \meter} & - \\
        $\phi$ & Produced/consumed energy & \unit{\watt} & $E$ \\
        $V_i$ & Control volume of node $v_i \in V$ & \unit{\cubic \meter} & $V$ \\
        \nrn{} & Number of nodes $\abs{V}$ & - & - \\
        \nre{} & Number of edges $\abs{E}$ & - & - \\
        $\mathcal{L}_e$ & Total energy/cost of operation & GWh/k\texteuro & -\\
        \hline
    \end{tabular}
\end{table}

As a precursor to the optimization problem, we first formulate our DHN in a mathematical sense.
The model that we use here to describe the DHN is similar to the one used in~\cite{krug_nonlinear_2021}.

\subsection{Network}
\begin{figure}[htb]
  \centering
  \includegraphics[width=\linewidth]{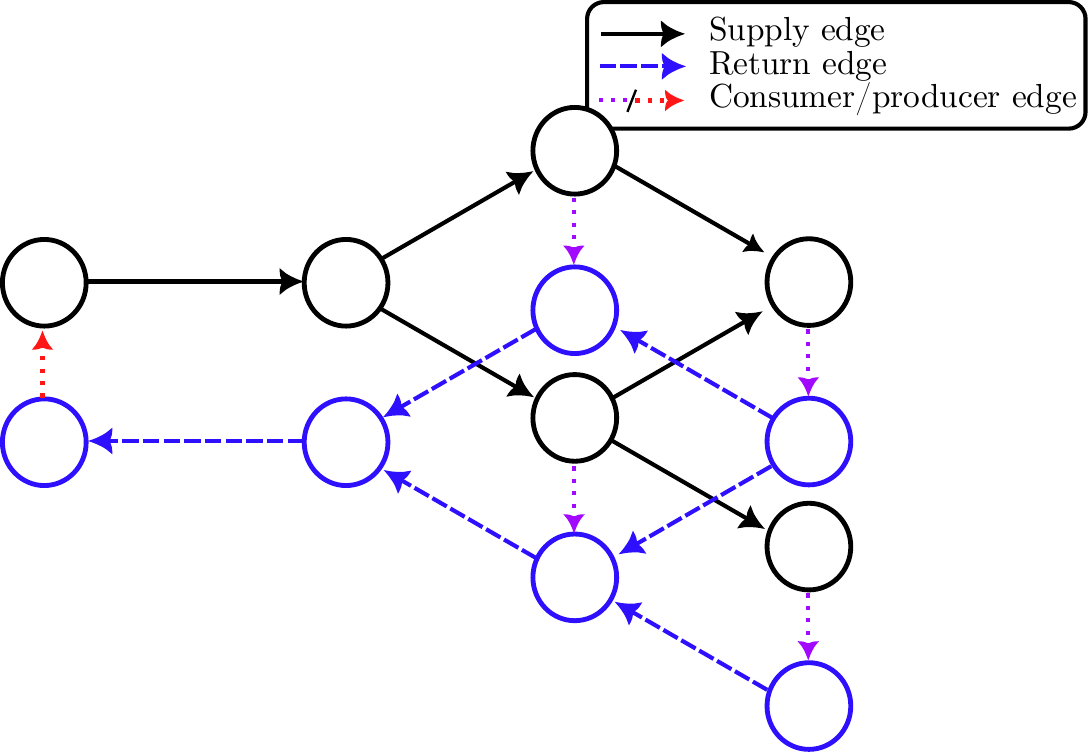}
  \caption{Illustration of a DHN graph, showing the duplicated graph structure and separation of nodes and edges.
  The nodes are separated into supply (black) and return nodes (blue), while the edges can either be considered supply (black), return (blue), consumer (purple) and producer (red) edges.
  Figure reproduced from~\cite{krug_nonlinear_2021}.
  }
  \label{fig:dhn_simple}
\end{figure}

The fundamental building block of modeling a DHN is by describing it as a connected, directed graph, consisting of vertices and edges, i.e.~$G = (V, E)$.
The graph is directed since we need to associate a flow direction across each edge.
Alternatively, $G$ can be expressed in terms of its incidence matrix $M_G \in \R^{\nrn{} \times \nre{}}$.
The graph actually consists of a duplicated structure, where we separate between the supply and return nodes, $V_s \cup V_r = V$ respectively.
Similarly, the edges $e = (v_i, v_j)$ for $v_i, v_j \in V$ are separated into supply $e_i \in E_s$ and return edges $e_i \in E_r$, connecting two supply or return nodes respectively.
Additionally however, we further have edges $e \in E_{c}$ connecting supply node with a return node.
The set of direct adjacent neighbors of a node $v_j \in \mathcal{N}(v_i)$, consists of all nodes connected to node $v_i$ by an edge.
A figure of a small, exemplary DHN graph is shown in \Cref{fig:dhn_simple}.
Further information on graph theory can be found in~\cite{bondy_graph_2008}.

For optimization, we consider a fixed time horizon on which we optimize $\mathbb{T} = \left[0, T \right], \, T>0$.
We consider multiple physical quantities that are modeled in the system, which are associated to either the edges $E$, the nodes $V$, or global constants.
A summary of the physical quantities and constants of this work is shown~\Cref{tab:variables}.

\subsection{Thermodynamic Model}
\label{sec:thermo_model}
In DHNs, power is distributed from producers of thermal energy via insulated pipes using fluids, most commonly.
The most prominent carrier of this heat exchange in DHNs is hot, pressurized water~\cite{mazhar_state_2018} with a thermal capacity of $c_p$. %
The heat transfer over the domain can be described using the convection-diffusion equation on a domain~\cite[Eq.~(21)]{guelpa_thermo_fluid_2019}
\begin{equation}
\rho c_p \dot{y} + \left<\rho c_p \vv, \nabla y \right> = \divop(k y) + \phi.
\label{eq:heat_cont}
\end{equation}
By assuming negligible diffusion in the flow direction (convection dominated) and by defining the control volume of a node $V_i$~\cite{guelpa_thermo_fluid_2019}, we write the discretized form as 
\begin{equation}
  \begin{cases}
    \rho c_p V_i \dot{y_i} + c_p \sum_{j} \dot{m}_j y_j = k l \left(y_i - y_a \right) + \phi & \quad i \notin BC  \\
    y_i = \hat{y}_i & \quad BC \subset \vcp{}.
    \end{cases}
\end{equation}
for $j$ being the edges adjacent to node $i$.
Note that the mass flow has become a scalar, positive only iff the fluid on an edge $e = (v_i, v_j) \in E$ is flowing from node $i$ to $j$, or negative otherwise.
To avoid stagnation, it is assumed that $\dot{m} \neq 0$ on any edges.
If we define $V = \operatorname{diag}(V_i)$, $G$ as the operator $\sum_j \dot{m}_j$ and discretizing the temperature $\vy$ as a vector, we can discretize the system using a backward Euler scheme as
\begin{equation}
  \begin{split}
    &\rho c_p V \frac{\vy_t - \vy_{t-\Delta_t}}{\Delta_t} + c_p G \vy_t + S \left( \vy_t - y_a \right)  \Phi = 0 \\
    \Leftrightarrow &\left(\frac{\rho c_p}{\Delta_t} V + c_p G + S \right) \vy_t = c \rho \vy_{t-\Delta_t} + S y_a 
  \end{split}
  \label{eq:heat_eq_lin_op}
\end{equation}
with additional boundary conditions at the plants and temperature deltas at the consumer edges representing $\phi$ in~\eqref{eq:heat_cont}
\begin{equation}
  \begin{cases}
    y_i = \hat{y}_i \quad & v_i \in BC \subset \vcp{} \\
    y_j = y_i - y_{d_i} \quad & (v_i, v_j) = e \in V_c 
  \end{cases}
  \label{eq:bc}
\end{equation}

\subsection{Optimal Control Problem}
\label{sec:ocp}
A popularized way of solving inverse problems in PDEs and ODEs is optimal control~\cite{manzoni_optimal_2021}.
In such a setting, the target function is minimized, while being constrained with the imposed PDE.

\paragraph{Solution Operator}
\label{sec:sol_op}
However, it is known that the thermodynamic model is well-posed given admissible boundary conditions and thus we can define a solution operator $S: U \to Y$ that directly maps from a given control vector $\vu$ (plant input temperatures over time) to a state vector $\vy$.
This greatly simplifies the optimal control we have to solve~\cite[p.~50]{troltzsch_optimal_2010}.

\paragraph{Loss}
\label{sec:loss}
We define a simple loss which reflects our physical quantity we want to minimize, the energy or money used in operating the DHN:
\begin{equation}
  \mathcal{L}_e(\vy) := c_p \sum_i \int_T \dot{m}_i(t) \left(y_{p_{s_i}} - y_{p_{r_i}} \right) p(t, y_{p_{s_i}}, y_{p_{r_i}}),
  \label{eq:loss}
\end{equation}
where we simply use $p(t, y_i, y_o) = 1$ for minimizing the energy.
Contrarily, when optimizing under varying energy prices (i.e.~electricity) $\tilde{p}: T \to \mathbb{R}$, we define our loss weighting as
\begin{equation}
  p(t, y_i, y_o) = 
  \begin{cases}
    \tilde{p}(t) \alpha &\text{if } y_i \ge y_o \\
    \tilde{p}(t) \beta &\text{if } y_i < y_o 
  \end{cases}
  \label{eq:gain_loss_facts}
\end{equation}
where the first case considers heat generation cost and the second a possible penalization when regaining energy from excessive heat in the network, assuming $\alpha, \beta > 0$.

\paragraph{Regularization}
\label{sec:regularization}
While our problem in \Secref{sec:sol_op} might be non-unique, we may impose smoothness on the solution through regularization methods~\cite{benning_regularization_2018}.
For the time dependant problem, we want to impose soft restrictions on the solution received from the optimization.
A first-order Tikhonov variant penalizes the quadratic variation of our control variable, effectively reducing rapid changes in the control solution:
\begin{equation}
  \mathcal{R}_t(u) := %
  \sum_{t} \left(\frac{\vu(t) - \vu(t-\Delta_t)}{\Delta_t}\right)^2.
\end{equation}

\paragraph{Constraints}
\label{sec:constraints}
In optimal control problems, constraints can either be enforced on the control $u$, or the state $y$.
In general, constraints in the control space $\Uad \subset U$ can be easily handled using projections~\cite[Sec.~6.4]{manzoni_optimal_2021} of the control vector $\vu$.
For state constraints $\Yad \subset Y$, we rely on penalty methods which we shortly describe here~\cite{boyd_convex_2004}:
Let us consider for such a case the inequality constraints and a simple loss function $f$
\begin{equation}
  \min_\vu f(\vu) \st c_i(\vu) \le 0,
  \label{eq:inequality_constr}
\end{equation}
which we alternatively formulate using the penalty method
\begin{equation}
  \min_\vu f_p(\vu) := \min_\vu f(\vu) + \frac{\lambda_p}{2} \dot{m}(c_i(\vu)),
  \label{eq:penalty}
\end{equation}
for 
$%
\dot{m}(c_i(\vu)) = \max(0, c_i(\vu))^2.
$ %
Note that for $\lambda_p \to \infty$,~\eqref{eq:penalty} and~\eqref{eq:inequality_constr} are equivalent minimization problems~\cite{boyd_convex_2004}.
In practice, we iteratively optimize the sequence of problems
\begin{equation}
  f_{p_i}(\vu) = f_{10 \lambda_{p_{i-1}}} (\vu),
  \label{eq:penalty_discrete}
\end{equation}
i.e.~we increase the penalization parameter by a factor of $10$ until at least $\lambda_p > 10^6$.

\subsection{Implementation}
The software is completely written in Python\footnote{\url{https://www.python.org}}, heavily utilizing PyTorch~\cite{paszke_pytorch_2019} for automatic differentiation of the losses in~\Secref{sec:loss}.
The automatic differentiation of the sparse solution operator in~\eqref{eq:heat_eq_lin_op} was implemented and published on github for this work\footnote{\url{https://github.com/thomgrand/torch_spsolve}}.
For the optimization of the loss in~\eqref{eq:loss}, we use L-BFGS~\cite{byrd_limited_1995} on the penalty loss.

\section{Experiments}
To test our methodology introduced in the last section, we will introduce the network and dynamic consumer load profiles that form the basis of our simulation experiment.
\subsection{Setup}

\paragraph{Network}
\begin{figure}
  \centering
  \includegraphics[width=\linewidth]{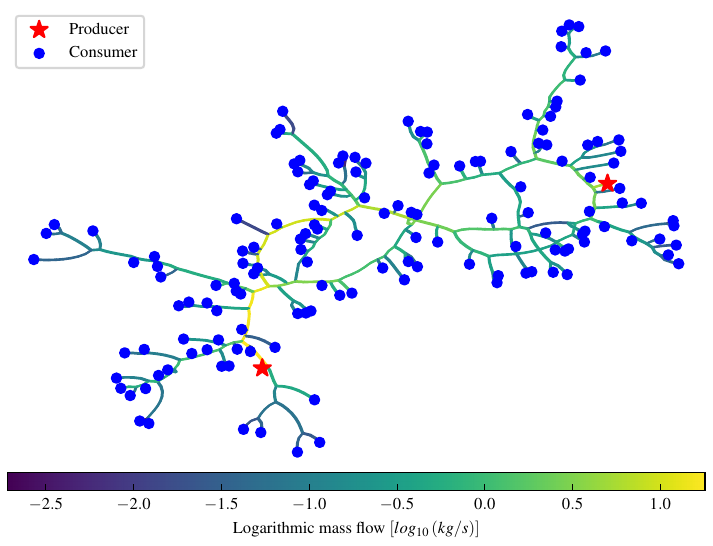}
  \caption{OpenDHN Net created from data of Verbier in Switzerland~\cite{boghetti_opendhn_2023}.
  The pipes of the DHN are shown as lines, the blue dots represent consumer sites and the red stars visualize the thermal supply plants.
  The color of the lines represents the mass flow through the network (logarithmic scale).
  }
  \label{fig:opendhn_net}
\end{figure}
As a basis for our simulation, we use an openly available DHN, called \emph{OpenDHN}~\cite{boghetti_opendhn_2023}.
This network is an abstraction of the Swiss alpine resort of Verbier, connected to 165 substations and 3 heating plants.
For privacy reasons, an equivalent layout of 150 substations and 2 heating plants is used.
This network consists of 1352 computational nodes, 1362 pipe segments, 150 consumer substations and 2 district heating plants.
For simulating the mass flows $\dot{m}$ throughout the network, given the consumer mass flows, we use the library \texttt{PyDHN}~\cite{boghetti_verification_2024}.
This steady mass flow is used throughout this work.
The network with the mass flows in logarithmic scale are shown in \Cref{fig:opendhn_net}.
It is important to note that we modified the mass flow rate of substation S132 to \SI{0.03}{\kilogram \per \second} in the network which was originally close to zero and also a source of major deviations in the original work~\cite[p.6]{boghetti_opendhn_2023}.

\begin{figure}[htb]
  \centering

  \includegraphics[width=\linewidth]{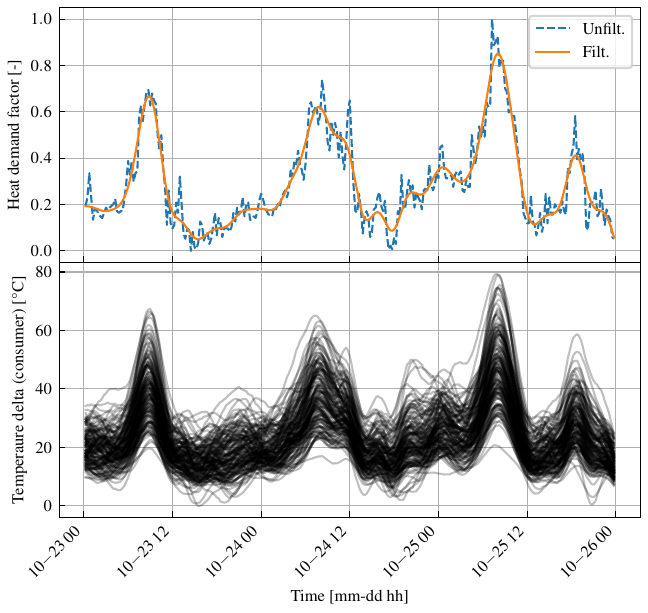}
  \caption{We show the average consumer model unfiltered in blue~\cite{de_mulder_2021_2022} and its low-pass filtered version in orange (top), along with the created individual variations of the single substations as transparent black lines (bottom), over the chosen three day period.
  }
  \label{fig:cons_model}
\end{figure}

\paragraph{Specific constraints}
\label{sec:spec_constraints}
Multiple constraints are considered in the optimization for differing reasons.
A direct constraint on the control is the maximum output temperature of each plant, which is limited to be below \SI{140}{\degreeCelsius}.
Especially in the dynamic energy pricing experiment, this will play a significant role, when the network is used as an energy storage in periods when energy is cheap.
Secondly, the supply temperature of the consumer needs to match a given minimal value of a given contractual value (\SI{80}{\degreeCelsius} in our case).
We additionally prescribe a minimum output temperature at each substation, since heat exchangers can only operate efficiently in certain temperature ranges, and our energy consumption at the consumer site is measured in terms of delta temperatures (see~\eqref{eq:bc}).
In short, the following the consumer input and output temperatures have to be fulfilled
\begin{equation}
  \begin{split}
    &y_{c, \text{s}} \ge \SI{80}{\degreeCelsius}  \\
    &y_{c, \text{r}} \ge \SI{30}{\degreeCelsius}  \\
    &y_{p} \le \SI{140}{\degreeCelsius}
  \end{split}
\end{equation}
\paragraph{Dynamic consumer model}
\label{sec:consumer_model}
On top of the provided network, OpenDHN offers also sensor date of the consumers and plants.
The data however has been anonymized, containing only averaged sensor data of mass flow, consumer energy demand and temperature supply.
In order to create a meaningful time dependant problem, we dynamically created synthetic consumer energy demand curves.
For this purpose, we started by utilizing a total heating district load profile over the period of three days of a Belgium ~\cite{de_mulder_2021_2022}.
As a period we chose the range from October 23 through October 25 of the year 2021.
This demand profile is further low-pass filtered with a 4th-order Butterworth with a cutoff frequency of roughly \SI{69.4}{\micro \hertz}.
From this average profile, we create multiple variations by varying the signal in Fourier space:
We vary the amplitude in a given frequency range by applying multiplicative noise to the Fourier domain of the timeseries before using an inverse Fourier transform to recover the varied timeseries signal.
In \Cref{fig:cons_model}, we visualize the total heat load loaded from~\cite{de_mulder_2021_2022}, the filtered profile and the inferred demand profile distribution of the 150 consumers.

\paragraph{Dynamic energy pricing}
\label{sec:dynamic_energy_pricing}
For our second experiment, we also consider dynamic energy pricing.
Such dynamic energy prices may arise from the integration of renewable energy sources into the thermal plant heat generation, directly optimizing economical viability of the DHN operation, or simply for penalizing operation during given periods or certain times of the day.

The dynamic energy prices are modeled by a simple time dependant function $\tilde{p}: T \to \mathbb{R}$, used in ~\eqref{eq:loss}.
Additionally, we distinguish the efficencies between a power-to-heat generation efficiency $\alpha$ and heat-to-power effiency $\beta$ through the function $p$~\eqref{eq:gain_loss_facts}.
In a case where no heat-to-power generation possibility is assumed, one can simply choose $\beta = 0$.

For the experiment at hand, we consider the hourly electricity day-ahead price of Switzerland at the European Power Exchange (EPEX Spot\footnote{\url{https://www.epexspot.com/en}}) for the chosen 3-day period~\cite{noauthor_european_nodate}, linearily interpolating in the 15 minute intervals we use for simulation and optimization.

\subsection{Results}
\begin{figure}[htb]
  \centering
  \includegraphics[width=\linewidth, trim={0cm 0cm 0cm 0cm}, clip]{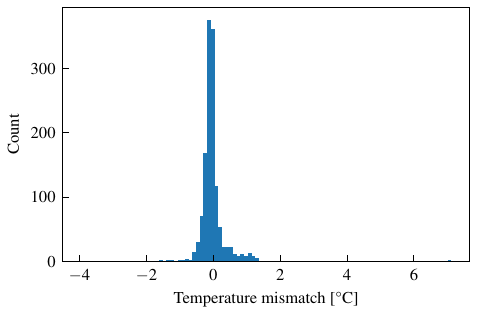}
  \caption{Histogram of the temperature mismatch ($x$-axis) of our solution operator $S$ against the \texttt{PyDHN} simulation on the OpenDHN benchmark.}
  \label{fig:pydhn_temp_mismatch}
\end{figure}
We first verified our thermal computations by testing the solution operator in steady state ($\dot{y}=0$) by applying the temperature deltas of each substations and plants according to the OpenDHN benchmark~\cite{boghetti_opendhn_2023} and computing the temperatures throughout the network.
The mean absolute mismatch in comparison to \texttt{PyDHN}~\cite{boghetti_verification_2024} was \SI{0.23}{\degreeCelsius}, with most values falling in the range \SIrange{-2}{2}{\degreeCelsius}.
A figure of the mismatch histogram is shown in \Cref{fig:pydhn_temp_mismatch}.

\begin{figure}[htb]
  \centering
  \includegraphics[width=\linewidth, trim={0cm 0cm 0.65cm 0cm}, clip]{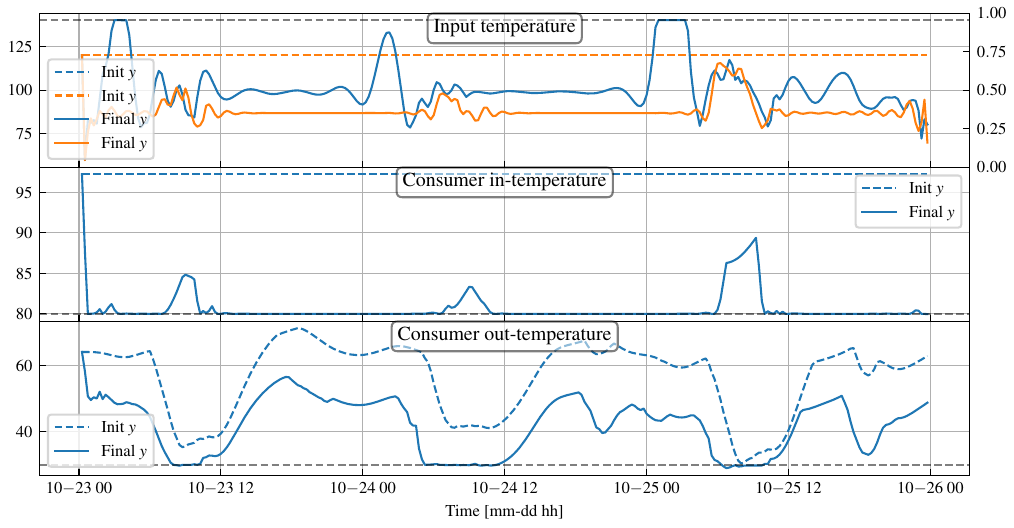}
  \caption{
    Optimal control result of the analyzed problem of Verbier over the chosen three day period.
    The top plot shows the input temperature of the network at the two thermal plants in different colors before (dashed) and after (solid) optimization.
    The middle plot shows the minimum input temperature over all consumer sites at the supply network.
    The bottom  plot shows the minimum output temperature over all consumer sites at the return network.
    All constraints shown are denoted in the plot as black dashed lines.
    All values in \unit{\degreeCelsius}.}
  \label{fig:results_in}
\end{figure}
\begin{figure}[htb]
  \centering
  \includegraphics[width=\linewidth, trim={0cm 0cm 0cm 0cm}, clip]{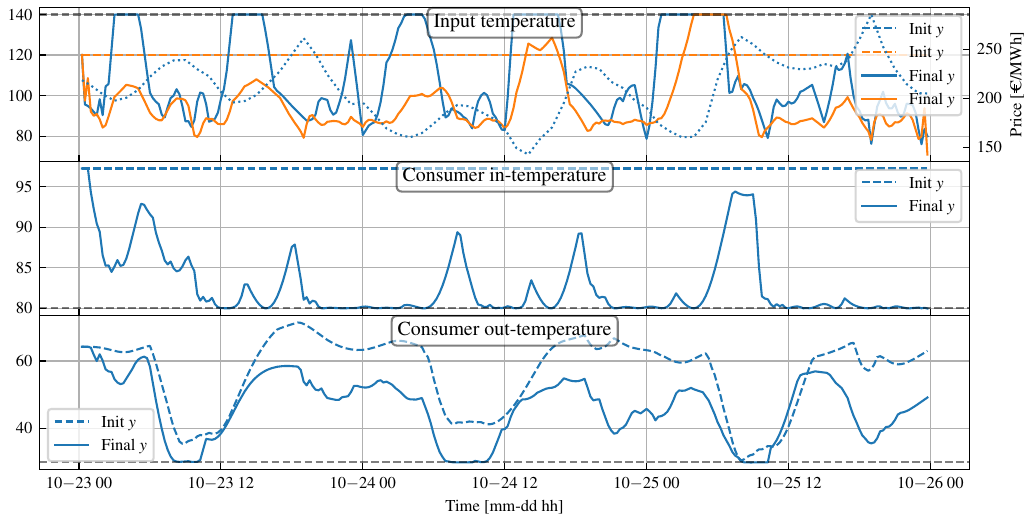}
  \caption{Same as \Cref{fig:results_in} for using our method to optimize for a time dynamic pricing model.
  The current energy price is shown additionally in the top plot (dotted blue line).}
  \label{fig:results_in_price}
\end{figure}
The results of our experiment of simulating and optimally controlling two plants in the OpenDHN heating grid are shown in \Cref{fig:results_in}.
Here, we show the temperatures at the two plants (top) before (dashed) and after (solid) optimization, as well as the minimal substation input and output temperature over all substations.
Since the energy consumed in the net is only the sum of ambient and exchange losses on top of the energy consumed by the substations, the natural optimal solution of \eqref{eq:loss} for a static energy price $p(t) = 1$ and given constraints is the low(est) temperature operation of the grid to minimize ambient losses.
As such, the input temperature is reduced to the minimal temperature that still fulfills the minimum substation input temperature or consumer output temperature, whichever is the restricting constraint.
In total, we can achieve an energy requirement reduction of 8\%.

\begin{figure}[htb]
  \centering
  \includegraphics[width=\linewidth, trim={0cm 0.85cm 0cm 0cm}, clip]{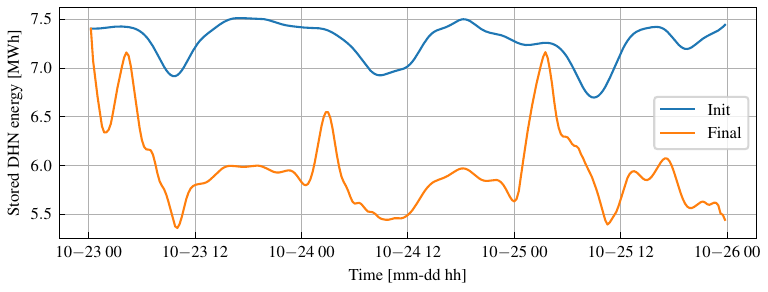}

  \vspace*{-0.075cm}

  \includegraphics[width=\linewidth, trim={0cm 0cm 0cm 0.075cm}, clip]{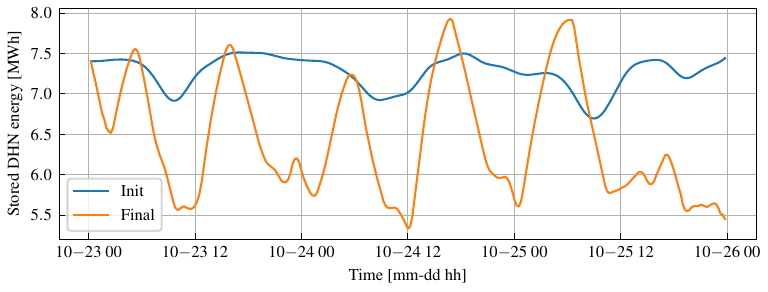}
  \caption{Stored energy inside the DHN for both tested scenarios over the chosen three day period.
  We show the stored energy over time for the low temperature operation optimal control (above) and for a dynamic pricing model (bottom).
  The blue line shows the initial operation and the orange line the optimized operation.}
  \label{fig:oc_energy_stored}
\end{figure}

By introducing dynamic energy pricing $p(t)$ in~\eqref{eq:loss} as described in \Secref{sec:dynamic_energy_pricing}, the optimal control drastically changes:
The optimization will put thermal energy into the network bounded by the maximum plant input temperature when energy is cheap, while lowering input temperatures at times of higher prices.
This can be seen by plotting the total amount of stored energy inside the DHN as was done in \Cref{fig:oc_energy_stored}.
The amount of stored energy exceeds the initial energy at times of low cost for cheaper consumer satisfaction.
It is notable, that towards the end of the three day period, the energy stored in the network will be depleted again in both scenarios, since charging the DHN at that point does not contribute to minimizing our loss.
In total, this optimization can achieve an energy cost reduction of 12.1\%.

\begin{figure}[htb]
  \centering
  \includegraphics[width=\linewidth]{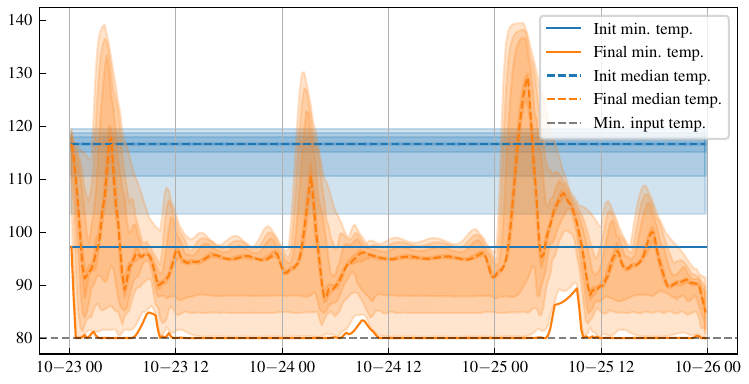}
  \caption{Quantile input temperature distribution at the consumer sites before (blue) and after (orange) optimization over the three days.
  The colored areas show the consumer supply temperatures for $\{10, 50, 90, 99\}\%$ of the consumers.
  The solid lines show the minimum and the dashed line the median supply consumer temperature.
  The grey dashed line marks the constrained minimum input temperature.}
  \label{fig:oc_temp_result_dist}
\end{figure}
\begin{figure}[htb]
  \centering
  \includegraphics[width=\linewidth]{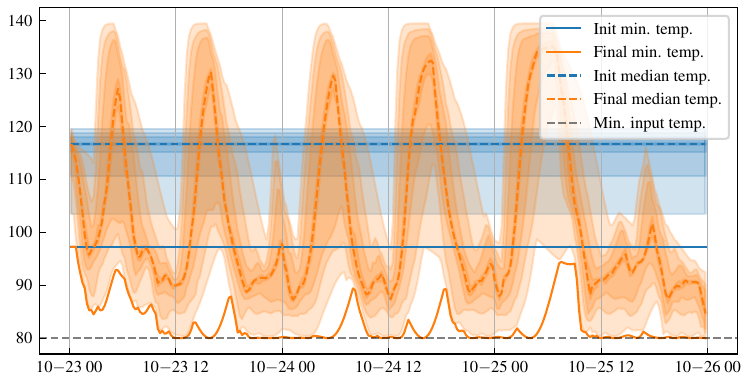}
  \caption{Same as \Cref{fig:oc_temp_result_dist} for the dynamic pricing model optimization.}
  \label{fig:oc_temp_result_dist_price}
\end{figure}
For further analysis, we look at the distribution of input temperatures at the substations.
In \Cref{fig:oc_temp_result_dist,fig:oc_temp_result_dist_price}, we show several quantiles (transparent area) and median (dashed color line) of the input temperatures at the substations before (blue) and after (orange) optimization.
The minimal input temperature constraint is shown as a dashed grey line, whereas the actual minimal input temperature is shown as a colored solid line.
We can see that especially in low temperature operation (\Cref{fig:oc_temp_result_dist}), the distribution of input temperatures is narrowed, with a difference of around \SI{15}{\degreeCelsius} between median and the 1\% temperature quantile and \SI{5}{\degreeCelsius} between median and the 99\% temperature quantile.
When using a dynamic pricing model, the input plant temperature varies a lot more as we have seen in \Cref{fig:results_in_price} and this is similarly reflected in a wider distribution of input temperatures at the substations.
In \Cref{fig:oc_temp_result_dist_price}, we see that every time the network is loaded with thermal energy during low energy prices, the difference in input temperatures at the substations increases up to more than \SI{50}{\degreeCelsius}.

\section{Discussion \& Conclusion}

We have shown that the optimization of a DHN can be performed using a solution operator combined with optimal control and its potential to optimize DHN thermal plant control to minimize either energy losses or associated costs.
The results show that the optimization reduces the energy consumption of the network by 8\% in a static energy price model and by 12.1\% in a dynamic energy price model.
The optimal control problem was shown in two different scenarios, optimizing a static or dynamic energy price model over a three-day period within a runtime of less than 5 minutes, making it feasible in real-time applications.
The utilized network, consumer models and dynamic energy pricing model were chosen and generated to reflect real-world scenarios as closely as possible, without the need to rely on sensitive, or not openly available data.

The consideration of dynamic pricing in DHNs reveals the potential to store excessively available energy in the form of thermal energy for a later demand surplus.
Using the excessively returned heat in such a system is not easily performed, as low temperature waste heat to power systems are often economically not viable.
However, they show great potential, especially in combination with thermal storages and might be a viable option in the future~\cite{xu_perspectives_2019}.

  \iftoggle{preprint}{
      \printbibliography %
  }{
      \bibliographystyle{ieeetr}
      \bibliography{thermal_demand_response}
  }

\end{document}